\documentclass{amsart}
\usepackage{amssymb}

\def\squareforqed{\hbox{\rlap{$\sqcap$}$\sqcup$}}
\def\qed{\ifmmode\squareforqed\else{\unskip\nobreak\hfil
\penalty50\hskip1em\null\nobreak\hfil\squareforqed
\parfillskip=0pt\finalhyphendemerits=0\endgraf}\fi\medskip}

\newcommand{\udot}{{}^{\textstyle .}}

\newcommand{\PSL}{\mathrm{L}}

\newcommand{\Mat}{\mathrm{M}}

\newcommand{\Sz}{\mathrm{Sz}}
\newcommand{\Co}{\mathrm{Co}}

\newcommand{\PSU}{\mathrm{U}}

\newcommand{\Ja}{\mathrm{J}}
\newcommand{\MM}{\mathbb{M}}

\title{On subgroups of the Monster isomorphic to $\PSL_2(8)$}
\date{Skeleton version begun 24th July 2016}
\author{Robert A. Wilson}

\address{School of Mathematical Sciences\\
Queen Mary University of London\\
London E1 4NS\\U.K.}
\email{R.A.Wilson@qmul.ac.uk}
\begin{document}
\maketitle

\begin{abstract}
We describe computer calculations that were used to classify
subgroups of the Monster isomorphic to $\PSL_2(8)$,
containing $7B$-elements. It turns out that there is no such
$\PSL_2(8)$ in the Monster.
\end{abstract}

\section{Introduction}
The maximal subgroup problem for 
the Fischer--Griess Monster group $\MM$ has been almost finished for several years. 
The papers \cite{oddlocals,MSh,Meier,Anatomy1,Anatomy2,post,A5subs,S4subs,L241,L227,L213B,Sz8A,U38}
have reduced the problem to classifying subgroups isomorphic to
one of the groups $\PSL_2(8)$, $\PSL_2(13)$, $\PSL_2(16)$ and $\PSU_3(4)$.
Of these, $\PSL_2(8)$ and $\PSL_2(16)$ have apparently
been classified in unpublished work
of P. E. Holmes, although the results do not appear
to be publicly available. 
Some information about the case $\PSL_2(13)$ is given elsewhere
\cite{L213A}.
In this paper we consider the case $\PSL_2(8)$, and give an
independent proof of the unpublished results of Holmes.
Throughout this paper, $\mathbb M$ denotes the Monster. Notation otherwise follows \cite{Atlas,FSG}, where much useful
background information can be found.

\section{Strategy}
Norton \cite{Anatomy1} has classified the subgroups $\PSL_2(8)$ 
of the Monster that contain
$7A$-elements. It remains therefore to classify those containing
$7B$-elements. Since $\PSL_2(8)$ contains a subgroup $2^3{:}7$, we first try to
locate suitable subgroups $2^3{:}7$ in the Monster.
Now the maximal $2$-local subgroups of the Monster have been
classified by Meierfrankenfeld and Shpectorov \cite{MSh,Meier}, and they are
just the ones listed in the Atlas
\cite{Atlas}. It is easy to see that the only one which
contains $7B$-elements is $2^{1+24}\udot \Co_1$. (This is 
essentially because a $7B$ element
does not centralize a $2A$-element.)

The elements of Monster class $7B$ in $2^{1+24}\udot\Co_1$ are of $\Co_1$ class $7A$.
To avoid confusion, we use the $\Co_1$ names for the classes here.
Moreover, the only maximal subgroups of $\Co_1$ which contain $7A$ elements
are the Suzuki chain subgroups, and hence there is no $2^3{:}7$ 
containing $7A$-elements in $\Co_1$. (A more elementary argument
is that any $2^3{:}7$ in $\Co_1$ lifts to $2\udot\Co_1$ and acts
faithfully on the Leech lattice, whence the element of order $7$
has fixed points so is in class $7B$.) In other words, 
every $2^3{:}7$ that we need to consider contains a
$2^3$ inside $2^{1+24}$.

Now the action of a fixed point free element of order $7$ on the $2^{24}$ is as the sum
of four copies of each of the $3$-dimensional representations. The $2^3$ we seek
lies in one of these homogeneous components, and it does not matter
which of the two homogeneous components
we take. The component is best thought of then as a $4$-space over the field of
order $8$, on which our element of order $7$ is acting as a scalar. The centralizer
$7\times A_7$ of the $7$-element acts irreducibly on this $4$-space, and it is an easy
calculation to show that it has three orbits on the $585$ subspaces of dimension $1$,
with lengths $15+210+360$. The stabilizers in $A_7$
are respectively $\PSL_3(2)$, $A_4$, $7$.

There are therefore just three groups $2^3{:}7$ that we need to consider.
We fix an element of order $7$, and compute representatives for the three
classes for $2^3{:}7$ that we need. We also compute the centralizer of this
element of order $7$, and an involution inverting it.
We are then able to find all the involutions inverting our element of order $7$,
and test each of the resulting amalgams of $2^3{:}7$ and $D_{14}$ to see whether
it generates a group $\PSL_2(8)$.
In more detail, the computation can be broken down into the following steps,
each of which is considered in a separate section below.
\begin{enumerate}
\item Finding $(A_7\times \PSL_2(7)){:}2$ in $\Co_1$.
\item Lifting to $(2\udot A_7\times \PSL_2(7)){:}2$ in $2^{1+24}\udot\Co_1$
\item Fixing an element of order $7$, and finding representatives for the
orbits of $7\times A_7$ on $1$-spaces (over the field of order $8$) in $2^{1+24}$.
We also need the stabiliser of each chosen representative.
\item Finding representatives for the orbits of each stabiliser on the 
involutions of cycle type $2^31$ in $S_7$.
\item `Changing post' to the centralizer of a new involution, that inverts
our chosen element of order $7$.
\item Identifying various elements already found with elements in the
new involution centralizer.
\item Searching in the new involution centralizer for elements which
centralize our chosen element of order $7$, in order to extend its centralizer
from the group $7\times 2\udot A_7$ already found to the full centralizer
$7^{1+4}{:}2\udot A_7$.
\item Computing a sufficient list of involutions inverting our element of
order $7$.
\item Checking the resulting list of cases to see if any of them is $\PSL_2(8)$.
\end{enumerate}
The paper is accompanied by a webpage 
where the necessary data and software are available:
{\tt http://maths.qmul.ac.uk/$\sim$raw/Monster/L28/}.
These can be used to check the calculations.
Notes throughout the text point to the relevant programs
for that section of text.

\section{Computing a copy of $(A_7\times \PSL_2(7)){:}2$ in $\Co_1$}
We look first for a dihedral group $D_{70}$, which necessarily contains
involutions from $\Co_1$ class $2C$.
To find a $2C$ element in $\Co_1$, we look for an element of order $22$.
We start from the `standard generators' $a$ and $b$, as defined in \cite{webatlas}.
Taking $c=ababab^2$ and $d=ab^2$ we find $e=cdcdcd^2cd$ has order $22$, so
$i_0=e^{11}$ is a $2C$-element. If $i_1=(i_0)^{ab}$, then $i_0i_1$ has order $35$,
and let $i_5=(i_0i_1)^{7}$ and $i_7=(i_0i_1)^{5}$.

Next we look for the centralizers of $i_5$and $i_7$, by searching for $3A$-elements
which commute with one or the other. The reason for adopting this strategy is that the
$3A$ conjugacy class is the smallest class in $\Co_1$, so this is
a reasonably efficient search.
To find a $3A$-element, look for an element of order $33$, such as
$f=abababab^2abab^2ab^2$, and then take $t_1=f^{11}$. We therefore look first for
a conjugate of $t_1$ which commutes with $i_5$. We find that $t_1$ conjugated by
$(ab^2)^7(ab)^{29}$, of order $66$, 
is such an element; call it $t_2$. Then $t_2$ and $i_7$ generate
$\Ja_2$, and $t_2,i_0,i_1$ generate $D_{10}\times \Ja_2$.
[Computational note: {\tt check1} computes the
elements {\tt c, d, e, i0, i1, i7, i5, f, t1, t2} from input {\tt a, b}.]

Next look for conjugates of $t_1$ which commute with $i_7$. We find $t_3$ and $t_5$,
defined as the conjugates of $t_1$ by $(ab)^{19}(ab^2)^{31}$, of order $28$, and by
$ab^2(ab)^{28}(ab^2)^{28}(ab)^{36}$, of order $24$, respectively. Each extends $i_5$ to $A_6$,
and together $t_3,t_5,i_5$ generate $A_7$, so that $t_3,t_5,i_0,i_1$ generate
$(A_7\times 7){:}2$.
[Computational note: {\tt check2} computes {\tt t3, t5} from {\tt a, b, t1}.]

To extend this to $(A_7\times \PSL_2(7)){:}2$, we look inside $\Ja_2$
for an appropriate subgroup $\PSL_2(7)$.
So we turn our attention to $\Ja_2$, and look first for a conjugate of $t_2$ commuting
with $t_3$. We find $t_6$, defined as $t_2$ conjugated by $(t_2i_7^2)^2$, of order $6$, and then
$t_6$ and $i_7$ generate $\PSU_3(3)$, commuting with the $A_6$ generated by $t_3$
and $i_5$. Then within $\PSU_3(3)$ we seek an involution commuting with $t_5$.
Such an involution is $t_7$, defined as the conjugate of $(t_6i_7)^2$ by
$(t_6i_7t_6i_7^2)^3$. We then have that $t_7,i_7$ generate $\PSL_2(7)$ commuting with the
$A_7$ generated by $t_3,t_5,i_5$.

Putting this together we have $(A_7\times \PSL_2(7)){:}2$ generated
by $i_0,i_1,t_3,t_5,t_7$. We also need an element of order $3$ normalizing $i_7$ in
the $\PSL_2(7)$ generated by $t_7,i_7$. Centralizing $i_0$, we find such an element
to be $t_8=(i_0t_7)^4(t_7i_7)^2$.

\section{Lifting to $2^{1+24}\Co_1$}
We now  lift the normalizer of $i_7$ to a group of shape $(2A_7\times 7{:}3).2$ in
$2^{1+24}\udot\Co_1$ by `applying the formula'. We lift the elements
$x=t_3$ and $t_5$, which centralize $i_7$ modulo the $2$-group, using the fact that
$(i_7i_7^x)^3$ conjugates $i_7^x$ to $i_7$, which gives the formula
$x'=i_7x(i_7i_7^x)^3$. (The initial $i_7$ is to ensure that $x$ and $x'$ lie in the same coset
of $2^{1+24}$.)
We adjust this formula for elements normalizing $i_7$ by
inserting the appropriate power, thus for $i_0$ we have
$i_0'=i_7i_0(i_7^{6}i_7^{i_0})^3$
and for $t_8$, which maps $i_7$ to its square, modulo the central
involution,
we have
$t_8'=i_7t_8(i_7^2i_7^{t_8})^3$.
(In fact, $i_7$ has order $14$ rather than $7$,
so we give the formulae we actually used, which may give 
different elements from the usual formulae.)
[Computational note:
{\tt check3} makes the elements
{\tt t6, t7, t8, t3p, t5p, i0p, t8p} from
the input {\tt t2, i7, i0, t3, t5}.]
\section{Finding the copies of $2^3{:}7$}
\label{find237}
The element $i_7^2$ of order $7$ acts on $2^{1+24}$ without fixed points. Its minimal
polynomial $x^6+x^5+x^4+x^3+x^2+x+1$ factorizes over the field of order $2$ as
$(x^3+x^2+1)(x^3+x+1)$. Its homogeneous components are found by taking the
nullspaces of these factors. The orbits of $7\times A_7$ are then easily found, and hence
we obtain vectors which give each of the three copies of $2^3{:}7$.
However, we want their stabilizers as well.

To find these, we make the permutation representation of $A_7$ on the 15 points
of the representation on the nullspace of $x^3+x+1$, and find the point stabilizer
generated by $t_9=(i_5t_3't_5'i_5t_3't_5'i_5)^3$ and $t_{10}=t_3't_5'i_5^2$.
Then the $2$-point stabilizer is generated by $t_9$ and $t_{11}=t_{10}^{6}t_9t_{10}^3$.
These two groups are the stabilizers of a $2^3$ in the first and second orbits
respectively. The third stabilizer is generated by $t_{10}i_7^2$.
Note, however, that by picking these particular stabilizers
we have made a choice of which of the two homogeneous components
we are working in. It turns out that the three cases we have
picked are not all in the same homogeneous component.

Now to make an element in one half of $2^{1+24}$, we take the element $a^2$,
which lies in $2^{1+24}$,
and compute $w=a^2i_7^{2}a^2i_7^{4}a^2i_7^{8}$;
while to get an element in the other half, compute
$w'=a^2i_7^{6}a^2i_7^{10}a^2i_7^{12}$. Then $i_8=(wt_{10})^7$
centralizes $t_{10}$, such that $i_7^2$and $i_8$
generate $2^3{:}7$ centralizing $\PSL_3(2)$.
Similarly, $i_{10}=(wt_{10}i_7^2)^7$ centralizes $t_{10}i_7^2$,
so that $i_{10}$ and $i_7^2$ generate $2^3{:}7$ centralizing
just $7$.

The remaining case proved slightly trickier to locate.
The vector $a^2$ turned out not to be good enough as a
seed vector, and we used instead the element $r_{10}$
defined in Section~\ref{7norm} below. Then let
$w''=r_{10}i_7^{6}r_{10}i_7^{10}r_{10}i_7^{12}$, and set
$i_9=(w''t_{11}^4)^3$. In fact,
we then need to replace $i_9$ by $i_9'=r_1^2i_9$,
so that it generates $2^3{:}7$
rather than $2^4{:}7$ with $i_7^2$. It is also necessary to check that
$i_9$ is not centralized by $t_{10}$ or its conjugates by 
$t_3^4t_5t_3^4t_5^2t_3^2$, to ensure that the $2^3$ is not centralized by a full
$2\udot \PSL_3(2)$, but only by $2\udot A_4$.

At the end of this computation, we have three non-conjugate copies of $2^3{:}7$,
being $\langle i_7,i_8\rangle$ and $\langle i_7,i_9\rangle$ and $\langle i_7,i_{10}\rangle$,
together with their centralizers. The centralizer of $\langle i_7,i_8\rangle$ is
$2\udot\PSL_3(2)$ generated by $t_{9}$
and $t_{10}$; the centralizer of $\langle i_7,i_9\rangle$ is
$2\udot A_4$ generated by $t_9$ and $t_{11}$; and the centralizer of $\langle i_7,i_{10}\rangle$ is
$14$, generated by $t_{10}i_7^2$ and the central involution.
[Computational note: {\tt check4}
makes the elements
{\tt t9, t10, t11, w, i8, i10} from input {\tt t3p, t5p, i5, a, i7}.
The element {\tt i9} is postponed until {\tt check11},
when {\tt r10} is available.]
\section{Finding orbit representatives on involutions}
There are $210$ conjugates of $i_0'$ in $2S_7$, and we need orbit representatives
for the action of each of the three stabilizers, that is $\PSL_3(2)$, $A_4$ and $7$
(modulo the central involution). Indeed, what we really want is elements of $S_7$ which
conjugate these orbit representatives to $i_0'$. This is a straightforward, if rather
tedious, computation in $2S_7$, so we merely present the results.

In the case of $\PSL_3(2)$, there are four orbits, of lengths $14$, $56$, $56$ and $84$, and
suitable elements to conjugate $i_0'$ to orbit representatives are
$$1, (i_0')^{t_3'}, t_3'^2, t_3'^2t_5'(t_3'^2)^{i_0'}.$$
Denote these four conjugating elements $o_0,o_1,o_2,o_3$
respectively.
In the case of $7$, there are $30$ orbits of length $7$. Six of them may be obtained
by conjugating first by $1$ or $(i_0')^{t_3'}$ and then by one of
$(t_3')^4, t_5',  t_3'^2t_5'(t_3'^2)^{i_0'}$. 
(These elements are denoted $o_{34}$ to $o_{39}$.)
The remaining $24$ are obtained by conjugating
first by $1$ or $(i_0')^{t'_3}$, then by one of $1,t_3'^2,t_5'^2,t_3'^2t_5'$, and
finally by one of $1,t_{11}t_9, (t_{11}t_9)^2$.
(These are denoted $o_{10}$ to $o_{33}$.)
In the case of $A_4$ we conjugate first by $1$ or $(i_0')^{t_3'}$, then by one of
$$1,(t_3')^4,(t_3')^2,t_5',(t_5')^2,(t'_3)^2t'_5, (t_3')^2t_5'((t_3')^2)^{i_0'}, t_{10},(t_3')^4t_{10},
(t_3')^2t_{10}$$
to get $20$ orbit representatives
(numbered from $o_{40}$ to $o_{59}$), and the remaining three from
$t_5't_{10}, (t_3')^4t_5't_{10}, (t_3')^2t_5't_{10}$
(numbered from $o_{60}$ to $o_{62}$).
[Computational note: {\tt check5}
computes these $57$ conjugating elements {\tt o0, ..., o62}
and their
inverses {\tt oi0, ..., oi62}.]

\section{Changing post}
For the next part of the calculation, we need to conjugate the inverting involution $i_0'$
to the central involution of $2^{1+24}\udot\Co_1$, so that we can find the rest of the
$7$-normalizer in the Monster.

 The first step is to work in the quotient $\Co_1$
to conjugate $i_0'$ into the normal $2^{11}$ of the standard copy
of $2^{11}{:}\Mat_{24}$. Now this standard copy is generated by $h$ and $i$ where
\begin{eqnarray*}
h&=&(ab)^{34}(abab^2)^3(ab)^6\cr
i&=&(ab^2)^{35}((ababab^2)^2ab)^4(ab^2)^5
\end{eqnarray*}
If we let
\begin{eqnarray*}
k_1&=&hihi^2\cr
k_2&=&hihihi^2\cr
k&=&(k_1k_2)^3k_2k_1k_2
\end{eqnarray*}
then $k$ has order
 $22$ in the quotient $\Co_1$, so its $11$th power maps to a
$2C$ element in $\Co_1$. 
We now find that 
$(k^{(ab)^6})^{11}i_0'$ has order $9$ modulo the central involution, and therefore
$$k_3=((k^{(ab)^6})^{11}i_0')^4(ab)^{-6}$$ conjugates $i_0'$ into the same
coset of $2^{1+24}$ that $k^{11}$ lies in.
A simple test then gives us that $T^{-1}$ conjugates 
this element $(i_0')^{k_3}$
into the normal $2$-group.

The second part of changing post is accomplished with the element
$$k_4=(ab)^{37}(ab^2)^{21}(ababab^2)^5(ab^2)^{28}(ab)^{10}.$$
Again, a simple test reveals we need to follow this with $T^{-1}$ rather than $T$.
That is, to conjugate $i_0'$ to the central involution $z$ we must conjugate by
$k_3T^{-1}k_4T^{-1}$.
[Computational note: 
{\tt check6} makes the elements
{\tt i, h, k1, k2, k, k3, k4} from input {\tt a, b, i0p}.]

\section{Moving to the new post}
We need to identify the old elements which centralize $i_0'$ with new elements in
the new post. The centralizer of $i_0'$ in $2A_7\times 7{:}3$ is $2A_4\times 3$,
generated by $t_8'$ and $s_1=i_0'(t_3'i_0')^3$ and 
$s_2=t_5't_3''$, where $t_3''$ is $t_3'$ conjugated by $(i_5t_3't_5'i_5)^3$,
of order $7$.
Writing $s_3=s_1t_8'$ gives us two generators $s_2$ and $s_3$.
We then write these as $24\times 24$ matrices over the field of order $3$, as
described in \cite{L241}, and work to identify them as elements
of the standard copy of $2\udot \Co_1$.

Since they are in the centralizer of $z'=z^{k_3T^{-1}k_4T^{-1}}$, they lie in
the standard copy of $2^{11}\Mat_{24}$, that is generated by $h$ and $i$.
We first identify the correct copy of $\Mat_{12}$ in $\Mat_{24}$ by finding two
elements $f_1,f_2$ of order $5$ which commute with the involution $z'$, such as
$(hi^2)^2$ conjugated by $(hi)^{10}(ih)^5$ and $(hi)^{18}(ih)^{11}$ respectively.
(In the Monster, these conjugating elements have orders $30$ and $88$
respectively.)
[Computational note: {\tt check7} makes the elements
{\tt s1, t3pp, s2, s3, f1, f2}
from input {\tt i0p, t3p, t5p, i5, t8p, h, i}.]

Then in the quotient $\Mat_{12}$ we identify the elements
$f_4=f_1f_2(f_1f_2f_1f_2^2)^2$ and $f_5=(f_1f_2)^4f_2f_1f_2$ of order $11$
to conjugate by, and then find that the image of $s_2^4$ may be obtained as
$(f_1f_2^2)^2$ conjugated by $f_4^8f_5^8f_4^9$, while the image of
$s_3$ may be obtained as $f_1f_2f_1f_2^2$ conjugated by
$f_4^5f_5^7f_4^8f_5^7$.

To lift back to $2^{11}\Mat_{12}$ we make generators $n_0=(f_1f_2)^8$
and $n_i=n_0^{f_4^i}$, and multiply on the left by
$n_0n_2n_4n_5n_6n_8n_9$ for $s_2^4$, and by
$n_3n_4n_5n_7n_8n_9$ for $s_3$. Finally we use our standard method to lift
back to $2^{1+24}\udot\Co_1$, multiplying on the left now by
$$d_2d_5d_6d_7d_9d_{10}d_{12}p_3p_4p_5p_6p_7p_9$$
for $s_2^4$, and by
$$d_1d_2d_3d_6d_7d_8d_9d_{12}p_1p_3p_5p_8p_9p_{10}p_{12}$$
for $s_3$.
Finally we check the results against the original elements, by applying the two versions
to a random vector, and find that $s_3$ is correct, and that $s_2^4$ needs further
adjustment by multiplying by the central involution $z$.
Denote the resulting elements by $u_2$ and $u_3$ respectively.
(In fact, the signs of the $p_i$ and the $d_i$ differ between different
versions, and so these overall signs are not reliable, and need to
be checked properly when required.)
[Computational note: 
{\tt check8} makes the elements
{\tt f4, f5, n0, ..., n9, u2, u3} from
the input {\tt f1, f2, d1, ..., d12, p1, ..., p12}.
The elements {\tt d1, ..., d12, p1, ..., p12} are made by {\tt check0}.]

\section{Finding the rest of the $7$-normalizer}
\label{7norm}
The next step is to extend the group $3\times 2A_4$ to $7^2{:}(3\times 2A_4)$ inside
the Conway group $\Co_1$. Now $3\times 2A_4$ has a central $3D$-element, $u_3^4$,
and
some $3B$ elements, including $u_2$, 
generating a particular copy of $2\udot A_4$. 
The remaining elements of order $3$ in $3\times 2A_4$
are in $\Co_1$ class $3D$.
To obtain $7^2{:}(3\times 2A_4)$, the plan is to find an element of order $7$
in the centralizer of one of these `diagonal' $3D$ elements, that is normalized by the
subgroup $3^2\times 2$. There are just four such $7$-cycles, up to symmetry,
two for each of the
diagonal elements. 

We first search for $3A$ elements in $\Co_1$ which centralize $u_2u_3^4$, modulo $2^{1+24}$.
We found $q_1$ and $q_2$, which are the conjugates of $t_1$ by
$(ab)^4(ab^2)^{17}(ab)^4$ and $(ab)^4(ab^2)^{32}(ab)^{38}$
respectively.
Then the group $3\times A_9$, centralizing $u_2u_3^4$, is generated by $u_3^2$ and $q_1q_2$.
By converting this group into the $9$ point permutation representation
of $A_9$, we can easily find
words for any desired element. The two generators can be mapped to the
permutations $(1,2,3,5,8,9)(4,6)$ and $(2,4,7)$ respectively. Then if $q_3,q_4,q_5$
are the conjugates of $q_1q_2$ by the square, cube and fifth power of $u_3^2q_1q_2$
respectively, we have the $7$-cycle $(q_3q_4q_5^2)^8$ mapping to $(5,7,3,8,2,1,9)$.
The conjugate of this by $q_5q_4q_5q_4^2q_5^2$ is then a $7$-cycle $q_6$
normalized by $u_3^2$. Conjugating $q_6$ by $q_4q_5q_4q_5^2$ gives $q_7$, the
other case.
[Computational note:
{\tt check9} makes the elements 
{\tt q1, q2, q3, q4, q5, q6, q7}
from the input {\tt a, b, t1, u3}.]

Now apply the formula to get elements which centralize $u_2u_3^4$
in $2^{1+24}\udot\Co_1$, that is
$q_6' = u_2u_3^4q_6u_2u_3^4(u_2u_3^4)^{q_6}$ and similarly for $q_7'$.
These lie inside the centralizer $2^{1+8}\udot A_9$ of $u_2u_3^4$.
Then apply another formula to get elements which are inverted by $u_3^6$, that is,
put $q_6''=q_6'^{-1}u_3^6q_6'u_3^6$ and similarly for $q_7''$.
These are not necessarily normalized by $u_3^4$, so we make further adjustments.
We calculate the centralizer of $u_3^6$ in the $2^{1+8}$ as follows.
Let $u_5=u_3^2$ and $u_6=q_1'q_2'$, where $q_1'$ and $q_2'$ are obtained by
applying the formula to $q_1$ and $q_2$, so they centralize $u_2u_3^4$.
Then $r_1=(u_5u_6u_5u_6^2)^9$ lies in the $2^{1+8}$, as do its conjugates
$r_{i+1}$ by $(u_5u_6)^i$, and the elements $r_i'=(r_iu_3^6)^2$ centralize $u_3^6$
modulo the central involution. In fact the centralizer of $u_3^6$
in $2^{1+8}$ is generated by the
central involution and $r_1'r_3', r_2',r_5'$.
The conjugates of $q_6''$ and $q_7''$ which are normalized by $u_5$
turn out to be $q_6''$ conjugated by $r_2'r_5'$ or $r_1'r_3'$, and $q_7''$
conjugated by $r_5'$ or $r_1'r_2'r_3'$. 

Finally we test these four elements of order $7$ to see if any of them commutes
with $i_7$ conjugated by $k_3T^{-1}k_4T^{-1}$. 
(In fact $i_7$ turned out to have order $14$, so we used $i_7^2$ instead.) 
It turned out that the last of the four cases worked, that is,
$q_7''$ conjugated by $r_1'r_2'r_3'$. Call this element $q_7'''$. 
At this stage we have generators for the full centralizer
of $i_7^2$ in the Monster. It is a group of shape 
$7^{1+4}{:}2A_7$ generated by $t_3', t_5', i_5$,
together with $q_7'''$ conjugated by $Tk_4^{-1}Tk_3^{-1}$.
[Computational note:
{\tt check10} makes the elements
{\tt q1p, q2p, q6p, q7p, q6pp, q7pp, u5, u6, r1, ..., r5,
r1p, ..., r5p, q7ppp}, from input {\tt u2, u3, q1, q2, q6, q7}.]

Of course, there was only a $50\%$ chance of success, and we actually calculated the other
case as well, that is, centralizing $u_2u_3^8$ rather than
$u_2u_3^4$. We found, of course 
as expected, that none of the elements obtained in
this case centralized $i_7^2$.
Specifically, we found centralizing elements
$q_{10}$ and $q_{11}$ being the conjugates of $t_1^4$ by
$(ab)^{10}(ab^2)^{26}(ab)^7$ and $(ab)^{18}(ab^2)^{21}(ab)^{17}$
respectively. Then $u_3^2$ and $q_{10}q_{11}$ generate the full $3\times A_9$ in
the quotient $\Co_1$. 
Applying the formula as above we obtain 
$q_{10}' = u_2u_3^8q_{10}u_2u_3^8(u_2u_3^8)^{q_{10}}$ 
and similarly for $q_{11}'$.
We define $u_7=q_{10}'q_{11}'$, and then $r_{10}
=((u_5u_7)^3u_7u_5u_7)^{21}$ is an involution which was used
(somewhat carelessly) as a seed for the element
$i_9$ defined in Section~\ref{find237}.
The rest of the calculation is suppressed, as it is not
required for the proof.
[Computational note:
{\tt check11} makes the elements
{\tt q10, q11, q10p, q11p, u7, r10, wpp, i9, i9p}
from input {\tt a, b, t1, u2, u3, u5, i7}.]

\section{Listing the involutions}
At this stage, we want to write down words for the $49$ conjugates of $a^4u_3^6$
under the $7^2$ group just constructed. However, the central involution in $3\times 2A_4$,
which acts as inversion on the $7^2$, also normalizes the groups $2^3{:}7$ constructed earlier,
so it is sufficient to take one involution in each orbit under this involution. Moreover, in the case
$a^4u_3^6$ itself, the group generated by $2^3{:}7$ and $7{:}2$ lies in
the involution centralizer, so this case can be ignored. 
Denote by $j_{7m+n}$ the conjugate of $j_0=a^4u_3^6$ by
$(q_7''')^n((q_7''')^{u_2})^m$. Then it is sufficient to use the cases $j_4,\ldots,t_{27}$.

(There is also in each case an element of order $3$ normalizing the group $2^3{:}7$, which 
could also be used to reduce the number of cases to consider. However, it is a different element
of order $3$ in each case, and it seemed safer not to use this symmetry.
In any case, the extra time taken in programming and debugging would likely outweigh the
saving.)

\section{Checking the cases}
In each case, we take the $7$ involutions in $2^3$, and conjugate by each of the
$57$ inverse conjugating elements in turn. Call the resulting elements $m_0,\ldots, m_6$.
Then we test (on a random vector) whether
$$Tk_4^{-1}Tk_3^{-1}mk_3T^{-1}k_4T^{-1}j$$
has order $3$, as $m$ runs through $m_0,\ldots, m_6$ and $j$ through $j_4,\ldots,j_{27}$.
It is a necessary and sufficient condition for the group generated to be $\PSL_2(8)$, that 
one of the $m_i$ gives a result of order $3$.

The total computation took around $36$ hours on a single processor on my (more than
ten years old) laptop. The final result was that no element of
order $3$ turned up. In other words, there is no $\PSL_2(8)$
containing $7B$-elements in the Monster.


\end{document}